\documentclass[11pt,a4paper]{article}
\pdfoutput=1 
\usepackage[british]{babel}
\usepackage{geometry}        
\usepackage[utf8]{inputenc}
\usepackage{graphicx} 
\usepackage{amssymb,latexsym,color,enumerate}
\usepackage[round]{natbib}
\title{\huge Some properties of skew-symmetric distributions}
\author{{\Large Adelchi Azzalini} \\
        Dipartimento di Scienze Statistiche\\
        Università di Padova\\
        Italia
        \and
        {\Large Giuliana Regoli}\\
        Dipartimento di Matematica e Informatica\\
        Università di Perugia \\
        Italia}
\date{21st December 2010}           

\usepackage{amsfonts}
\newcommand{\anymode}[1]{\ifmmode{#1}\else\mbox{$#1$}\fi}

\newcommand{\NB}[1]{\texttt{/\textsuperscript{NB}}%
          \marginpar{\raggedright\sl\small#1\hfill}}	  

\long\def\Nota#1{\footnote{#1}\kern-0.2em\NB{cfr Nota$^{(\thefootnote)}$}}
\def\linea#1{\ifhmode\hfill\break\fi\hbox to \hsize{#1}}

\newcommand{\inv}{^{-1}}

\def\vc{v\kern -0.1em .c.\relax}

\newcommand{\dfrac}[2]{\displaystyle{\frac{#1}{#2}}}

\newcommand{\equald}{\stackrel{d}{=}}
\newcommand{\E}[2][]{
   \ensuremath{\mathbb{E}_{#1}\!\left\{\displaystyle{#2}\right\}}}

\newcommand{\half}{\mbox{$\textstyle \frac{1}{2}$}}   
\long\def\ignore#1{}

\newcommand{\indep}{\perp\kern-0.5em\perp}

\newcommand{\pr}[2][]{
   \ensuremath{\mathbb{P}_{#1}\!\left\{\displaystyle{#2}\right\}}}

\newcommand{\Real}{\mathbb{R}}

\newcommand{\sign}{\mathop{\rm sgn}\nolimits}

\newcommand{\T}{^{\top}}
\newcommand{\var}[2][]{
   \ensuremath{\textrm{var}_{#1}\!\left\{\displaystyle{#2}\right\}}}

\newtheorem{theorem}{Theorem}
\newtheorem{lemma}[theorem]{Lemma}
\newtheorem{definition}[theorem]{Definition}

\newtheorem{corollary}[theorem]{Corollary}
\newtheorem{proposition}[theorem]{Proposition}

\renewcommand{\d}{\,\mathrm{d}}
\newcommand{\pd}[2]{\displaystyle\frac{\partial #1}{\partial #2}}

\newcommand{\GR}{\ensuremath{\ge_{\hbox{\tiny GR}}}}
\begin{document}
\maketitle
\thispagestyle{empty}
\begin{abstract}
  The family of skew-symmetric distributions is a wide set of probability
  density functions obtained by combining in a suitable form a few components
  which are selectable quite freely provided some simple requirements are
  satisfied.  Intense recent work has produced several results for specific
  sub-families of this construction, but much less is known in general terms.
  The present paper explores some questions within this framework, and
  provides conditions on the above-mentioned components to  ensure that the 
  final distribution enjoys specific properties.
\end{abstract}

\noindent\emph{Primary Subjects}: 60E05 62H10

\noindent\emph{Secondary Subjects}: 62E15, 60E15


\clearpage
\section{Introduction and motivation}
\subsection{Distributions generated by perturbation of symmetry}
In recent years, there has been quite an intense activity connected to a broad
class of continuous probability distributions which are generated starting
from a symmetric density functions and applying a suitable form of
perturbation of the symmetry.  The key representative of this formulation is
the so-called skew-normal distribution, whose density function in the
scalar case is given by
\begin{equation} \label{e:SN-pdf}  
    f(x;\alpha) = 2\,\phi(x)\,\Phi(\alpha\,x),  \qquad\qquad (x\in\Real),
\end{equation}
where $\phi(x)$ and $\Phi(x)$ denote the N(0,1) density function and
distribution function, respectively, and $\alpha$ is an arbitrary real
parameter.  When $\alpha=0$, (\ref{e:SN-pdf}) reduces to the N(0,1)
distribution; otherwise an asymmetric distribution is obtained, with skewness
having the same sign of $\alpha$.  Properties of (\ref{e:SN-pdf}) studied by
\citet{azzalini:1985} and by other authors show a number of similarities with
the normal distribution, and support the adoption of the name skew-normal.

Furthermore, the same sort of mechanism leading from
the normal density function to (\ref{e:SN-pdf}) has been applied to other
symmetric distributions, including extensions to  more elaborate forms 
of perturbation and constructions in the multivariate setting.
Introductory accounts to this research area are provided by the book edited 
by \citet{genton:2004-SE} and by the review paper of \citet{azzalini:2005}, 
to which the reader is referred for a general overview.

For the aims of the present paper, we shall largely rely on the following
lemma, presented by \citet{azza:capi:2003}.  This is very similar to an
analogous result developed independently by \citet{wangJ:boye:gent:2004-ss};
the precise interconnections between the two statements will be discussed in
the course of the paper. Before stating the result, we recall that the notion
of symmetric density function has a simple unique definition only in the
univariate case, but in the multivariate case there exist different
formulations; see \citet{serfling:2006-ess2} for an overview.  In this paper,
we adopt the notion of central symmetry, which in the case of a continuous
distribution on $\Real^d$ requires that a density function $p$ satisfies
$p(x-x_0)=p(x_0-x)$ for all $x\in\Real^d$, for some centre of symmetry $x_0$.

\begin{lemma} \label{th:prop1} 
Denote by $f_0$ a $d$-dimensional probability density function centrally
symmetric about 0, by $G_0(\cdot)$ a continuous
distribution function on the real line such that $g_0=G_0'$ is an even
density function, and by $w$ an odd real-valued function on $\Real^d$ 
such that $w(-x)=-w(x)$.  Then 
\begin{equation} \label{e:prop1}
       f(x) = 2\,f_0(x)\,G_0\{w(x)\}, \qquad\qquad (x\in\Real^d),
\end{equation}  
is a density function.
\end{lemma}

This result provides a general mechanism for modifying an initial symmetric
`base' density $f_0$ via the perturbation factor $G(x)=G_0\{w(x)\}$, whose
components $G_0$ and $w$ can be chosen among a wide set of options.  Clearly,
the prominent case (\ref{e:SN-pdf}) can be obtained by setting $d=1$,
$f_0=\phi$, $G_0=\Phi$, $w(x)=\alpha\,x$ in (\ref{e:prop1}).
The term `skew-symmetric' is often adopted for distributions of type
(\ref{e:prop1}).  An important property associated to Lemma~\ref{th:prop1} is
provided by the next statement.

\begin{proposition}[Perturbation invariance] \label{th:prop2}
If the random variable $X_0$ has density $f_0$ and $X$ has density $f$,
where $f_0$ and $f$ satisfy the conditions required in Lemma~\ref{th:prop1}, 
then the equality
\begin{equation} \label{e:prop2}
   t(X)  \equald t(X_0) \,,
\end{equation}  
where `$\equald$' denotes equality in distribution, holds
for any even $q$-dimensional function $t$ on $\Real^d$, irrespectively
of the factor $G(x)=G_0\{w(x)\}$.
\end{proposition}
 

\subsection{A wealth of open questions}
 
The intense research work devoted to distributions of type (\ref{e:prop1}) 
has provided us with a wealth of important results. 
Many of these have however been established for specific subclasses of 
(\ref{e:prop1}). The most intensively studied instance is given by the 
skew normal density which in the case $d=1$ takes the form (\ref{e:SN-pdf}). 
Important results have been obtained also for other subclasses,
especially when $f_0$ is the Student's $t$ density or the Subbotin 
density (also called exponential power distribution).

Much less is known in general terms, in the sense that there still is a
relatively limited set of results which allow us to establish in advance, on
the basis of qualitative properties of the components $f_0, G_0, w$ of
(\ref{e:prop1}), what will be the formal properties of the resulting density
function $f$.  Results of this kind do exist, and Proposition~\ref{th:prop2}
is the most prominent example, since it is both completely general and of
paramount importance in the associated distribution theory; from this
property, several results on quadratic forms and even order moments follow.
Little is known about the distribution of non-even transformations.
Among the limited results of the latter type, some general properties of odd
moments of (\ref{e:prop1}) have been presented by
\citet{umbach:2006,umbach:2008}.  There are however many other questions,
which arise quite naturally in connection with Lemma~\ref{th:prop1}; the
following is a non-exhaustive list.
\begin{itemize}
\item
  In the case $d=1$, which assumptions on $G(x)$ ensure that the median 
  of $f$ is larger than 0? More generally, when can we say the the
  $p$-th quantile of $f$ is larger than the $p$-th quantile of $f_0$?
  Obviously, `larger' here can be replaced by `smaller'.
\item
  The even moments of $f$ and those of $f_0$ coincide, because of 
  (\ref{e:prop2}). What can be said about the odd moments?
  For instance, is there an ordering of moments associated to some
  form of ordering of $G(x)$?
\item 
  If $f_0$ is unimodal, which are the additional assumptions on
  $G_0$ and $w$ which ensure that $f$ is still unimodal? 

\item When $d>1$, a related but distinct question is whether high density
  regions of the type $C_u=\{x : f(x)>u\}$, for an arbitrary positive $u$,
  are convex regions.
\end{itemize}

The aim of the present paper is partly to tackle
the above questions, but at the same time we take a broader view, attempting
to make a step forward in understanding the general properties of the set
of distributions (\ref{e:prop1}).  The latter target is the motivation for the
preliminary results of Section\,\ref{s:common-base}, which lead to a
characterization result in Section\,\ref{s:char} and provide the basis for the
subsequent sections which deal with more specific results.  In
Section~\ref{s:d=1} we deal with the case $d=1$ and tackle some of the
questions listed above.  Specifically, we obtain quite general results on
stochastic ordering of skew-symmetric distributions with common base $f_0$,
and these imply orderings of quantiles and of expected values of suitable
transformations of the original variate.  The final part of
Section~\ref{s:d=1} concerns uniqueness of the mode of the density $f$.
Section~\ref{s:d-general} deals with the case of general $d$, where various
results are obtained. One of these is to establish convexity of the sets $C_u$
for the more important subclass of the skew-elliptical family, provided the
parent elliptical family enjoys the same property. We also examine the
connection between the formulation of skew-elliptical densities of type
(\ref{e:prop1}) and those of \citet{bran:dey:2001}, and prove the conjecture
of \citet{azza:capi:2003} that the first formulation strictly includes the
second one. Finally we gives conditions for the log-concavity of
skew-elliptical distributions not generated by the conditioning mechanism of
\citet{bran:dey:2001}.

\section{Skew-symmetric densities with a common base}
\label{s:common-base}

\subsection{Preliminary facts}   
\label{s:equivalence}

Clearly, $f$ in (\ref{e:prop1}) depends on $G_0$  only via
the perturbation function $G(x)=G_0\{w(x)\}$. The assumptions on
$G_0$ and $w$ in Lemma~\ref{th:prop1} ensure that 
\begin{equation}  \label{e:G}
           G(x)\ge 0 \,, \qquad G(x)+G(-x) = 1 \,,\qquad(x\in\Real^d),
\end{equation}
and it is conversely true that a function $G$ satisfying these
conditions ensures that
\begin{equation}  \label{e:prop1-G}
           f(x) = 2\,f_0(x)\,G(x) 
\end{equation}
is a density function. In fact (\ref{e:G})--(\ref{e:prop1-G}) represent
the formulation adopted by \citet{wangJ:boye:gent:2004-ss} for their 
result essentially equivalent to Lemma~\ref{th:prop1}.

Each of the two formulations has its own advantages. As remarked by 
\citet{wangJ:boye:gent:2004-ss}, the representation of $G(x)$ in
the form $G(x)=G_0\{w(x)\}$ is not unique. In fact, given one
such representation,  
\[   G(x) = G_*\{w_*(x)\}, \quad w_*(x)= G_*\inv[G_0\{w(x)\}] \]
is another one, for any strictly increasing distribution function $G_*$
with even density function on $\Real$. 

On the other hand, finding a function $G$ fulfilling conditions (\ref{e:G}) is
immediate if one builds it via the expression $G(x)=G_0\{w(x)\}$; in fact, this
is the usual way adopted in the literature to select suitable $G$ functions. 
Furthermore, \citet{wangJ:boye:gent:2004-ss} have shown that the converse 
fact holds: any function $G$ satisfying (\ref{e:G}) can be written
in the form $G_0\{w(x)\}$,  and this can be done in infinitely many ways. 
A choice of this representation which we find `of minimal modification' is 
\begin{equation}  \label{e:G0-G}
   \begin{array}{lc}
    G_0(t)=\left(t+\half\right)\,I_{(-1,1)}(2t)+I_{[1,+ \infty )}(2t)\,, & 
        (t\in \Real)\,, \\[1ex]
   w(x)= G(x)-\half \,, & ( x\in \Real^d)\,,
   \end{array}
\end{equation} 
where $I_A(x)$ denotes the indicator function of the set $A$. In plain words,
this $G_0$ is the distribution function of a $U(-\half,\half)$ variate.

Another important finding of \citet[Proposition 3]{wangJ:boye:gent:2004-ss} 
is that  any positive density function $f$ on $\Real^d$ admits a 
representation  of type (\ref{e:prop1-G}), as indicated in their result
which we reproduce next with a little modification concerning the
arbitrariness of $G(x)$ outside the support of $f_0$.
Here and in the following, we denote by $-A$ the set formed by reversing the 
sign of all elements of $A$, if $A$ denotes a subset of a Euclidean space.
If $A=-A$, we say that $A$ is a symmetric set.

\begin{proposition} \label{th:propW} Let $f$ be a density function with
  support $S\subseteq\Real^d$. Then a representation of type (\ref{e:prop1-G})
  holds, with
  \begin{equation} \label{e:Wang}
  \begin{array}{rcl}
   f_0(x) &=&  \cases{\half\{f(x)+f(-x)\} & if $x\in S_0$, \cr
                               0 & otherwise,} \\[2.5ex]
    G(x)  &=&  \cases{\dfrac{f(x)}{2 f_0(x)} & if $x\in S_0$, \cr
                      \mathrm{arbitrary} & otherwise,}                
  \end{array}
\end{equation}
where $S_0=(-S)\cup S$, and the arbitrary branch of $G$ satisfies (\ref{e:G}).
Moreover $f_0$ is unique, and $G$ is uniquely defined over $S_0$.
\end{proposition}

Consider now a density function with representation of type (\ref{e:prop1-G}).
We first introduce a property of the cumulative distribution function 
$F$  which is also of independent interest.
Rewrite the first relation in (\ref{e:Wang}) as
\begin{eqnarray} 
     f(-x) = 2\,f_0(x) - f(x). \label{e:f(-x)}
\end{eqnarray}
for any $x=(x_1,\dots,x_d)$.
 If we denote by $F_0$ the cumulative distribution function of $f_0$, 
then integration of (\ref{e:f(-x)}) on  $\cap_{j=1}^d (-\infty,x_j]$ gives
\begin{equation} \label{e:F(x.)}
   \overline{F}(-x) =   2\,F_0(x) - F(x)   \,
\end{equation}
where $\overline{F}$ denotes the survival function, that is 
$\overline{F}(x)=\pr{X_1 \ge x_1, \dots, X_d\ge x_d}$;
(\ref{e:F(x.)}) can be written as
\[ \overline{F}(-x) + F(x) = F_0(x) +\overline{F_0}(-x) \] 
and  
this is in turn equivalent to Proposition~\ref{th:prop2}, 
as stated in Proposition~\ref{th:char} below.

\subsection{A characterization} \label{s:char}

The five single statements composing the next proposition are known
for the case $d=1$, some of them also for general $d$. The more 
important novel fact is their equivalence, which therefore represents
a characterization type of result.
 
\begin{proposition}\label{th:char}
  Consider a  random variable $X=(X_1,...X_d)\T$ with density
  function $f$ and cumulative distribution function $F$, and a continuous
  random variable $Y=(Y_1,...Y_d)\T$ with density function $h$ and 
  distribution function $H$. Then the following conditions are equivalent:
\begin{enumerate}[\quad (a)]
\item
  the densities $f(x)$ and $h(x)$ admit a representation of type
  (\ref{e:prop1-G}) with the same symmetric base density $f_0(x)$,
  
\item  $t(X)\equald t(Y)$, 
  for any even $q$-dimensional function $t$ on $\Real^d$, 

\item $P(X\in A) = P(Y\in A)$, for any symmetric set $A \subset \Real^d$,

\item $F(x) + \overline{F}(-x) =H(x) + \overline{H}(-x)$,

\item $ f(x)+f(-x) = h(x)+h(-x), \quad\hbox{(a.e.)}$.
\end{enumerate}
\end{proposition}

\noindent\textbf{Proof}
\begin{description}
\item[(a)$\Rightarrow$(b)] This follows from
the perturbation invariance property  of Proposition~\ref{th:prop2}.

\item[(b)$\Rightarrow$(c)] Simply notice that  the indicator function 
of a symmetric  set $A$  is an even function.

\item[(c)$\Rightarrow$(d)] 
On setting
\begin{eqnarray*}
   A_{+}&=& \{s=(s_1,\dots,s_d)\in \Real^d : 
   s_j\leq x_j, \forall j \}, 
    \\
   A_{-}   &=&  \{s=(s_1,\dots,s_d)\in \Real^d : 
   -s_j\leq x_j, \forall j \} 
      =  -A_{+}, \\
  A_\cup &=& A_{+}\cup A_{-} \,,\\
  A_\cap &=& A_{+}\cap A_{-}\,,
\end{eqnarray*}
both $A_\cup$,  $A_\cap$ are symmetric sets; hence we get
\begin{eqnarray*}
   F(x) + \overline{F}(-x)  
      &=& P(X\in A_{+}) +P(X\in A_{-}) \,, \\
      &=& P(X\in A_\cup) + P(X\in A_\cap) \,.
\end{eqnarray*}      

\item[(d)$\Rightarrow$(e)] 
Taking the $d$-th mixed derivative of (d), relationship (e) follows.

\item[(e)$\Rightarrow$(a)] It follows from the representation given in 
Proposition \ref{th:propW}.
\end{description}

In the special case $d=1$, the above statements can be re-written in more 
directly interpretable expressions.  Specifically, (\ref{e:F(x.)}) leads to
\begin{equation} \label{e:F(-x.)}
   1-F (-x) =   2\,F_0(x) - F(x)   \,,
\end{equation}
which will turn out to be useful later, and
\[  F(x) - F(-x) = F_0(x) - F_0(-x). \]

Moreover, when $d=1$, conditions (c) and (d) in Proposition \ref{th:char} 
can be replaced  by the following more directly interpretable forms:
\begin{description}
\item[\quad (c$'$)] $ |X| \equald|Y|$,
\item[\quad (d$'$)]
$  F(x) - F(-x) = H(x) - H(-x)$\,,
\end{description}
the first of which has appeared in \citet{azzalini:1986}, and the second
one is an immediate consequence.

\section{Some results when $d=1$}  \label{s:d=1}

\subsection{Stochastic ordering the univariate case}

In this section, we focus on the case with $d=1$.  We first introduce an
ordering on the set of functions which satisfy (\ref{e:G}). When this concept
is restricted to symmetric distribution functions, it reduces to the peakedness
order introduced by \citet{birnbaum:1948}, to compare the variability of
distributions about $0$.

\begin{definition} 
  If $G_1$ and $G_2$ satisfy (\ref{e:G}), we say that $G_2$ is greater than 
  $G_1$ on the right, denoted $G_2 \GR G_1$, if $G_2(x)\ge G_1(x)$ for all 
  $x>0$ and  strict inequality holds for some $x$.
\end{definition}
Of course it is equivalent to require that $G_2(x)\le G_1(x)$ for all $x<0$ 
and the inequality holds at some $x$.  Another equivalent condition is that
\[   G_2(s)-G_2(r)\ge G_1(s)- G_1(r), \qquad(r<0<s)\,. \]

If we now consider a fixed symmetric `base' density $f_0$ and the perturbed 
distribution functions associated to $G_1$ and $G_2$, that is
\begin{equation} \label{e:F}
    F_k(x) = \int_{-\infty}^x 2\,f_0(u)\, G_k(u) \d{u}, \qquad (k=1,2),
\end{equation}
the ordering $G_2 \GR G_1$ implies immediately the stochastic ordering of
$F_1$ and $F_2$ in the usual sense that $F_2$ is stochastically larger than
$F_1$ if $F_1(s)\ge F_2(s)$ for all $s$.  To see this, consider first $s\leq
0$; then $G_1(x)\ge G_2(x)$ for all $x\le s$, and this clearly implies
$F_1(s)\ge F_2(s)$. If $s>0$, the same conclusion holds by using
(\ref{e:F(-x.)}) with $x=-s$. We have then reached the following conclusion.

\begin{proposition} 
\label{th:ordering}
If $G_1$ and $G_2$  satisfy condition (\ref{e:G}), and
$G_2 \GR G_1$, then the distribution functions (\ref{e:F}) satisfy
\begin{equation} \label{e:ordering}
  F_1(x) \ge F_2(x) \,,\qquad (x\in\Real).
\end{equation}
\end{proposition}

Since $G_0$ is a monotonically increasing function, then it can be easier
to check the ordering of $G_1$ and $G_2$ via the ordering of the
corresponding $w(x)$'s. 

\begin{proposition} 
\label{th:ordering-w}
If $G_1=G_0(w_1(x))$ and $G_2=G_0(w_2(x))$  where $G_0$ is as in 
Lemma~\ref{th:prop1},  and $w_1$ and $w_2$ are odd functions 
such that $w_2(x)\ge w_1(x)$ for all $x>0$, then $G_2 \GR G_1$
and (\ref{e:ordering}) holds. 
\end{proposition}
 
\begin{figure}[htbp] 
   \centering
   \includegraphics[width=0.48\hsize]{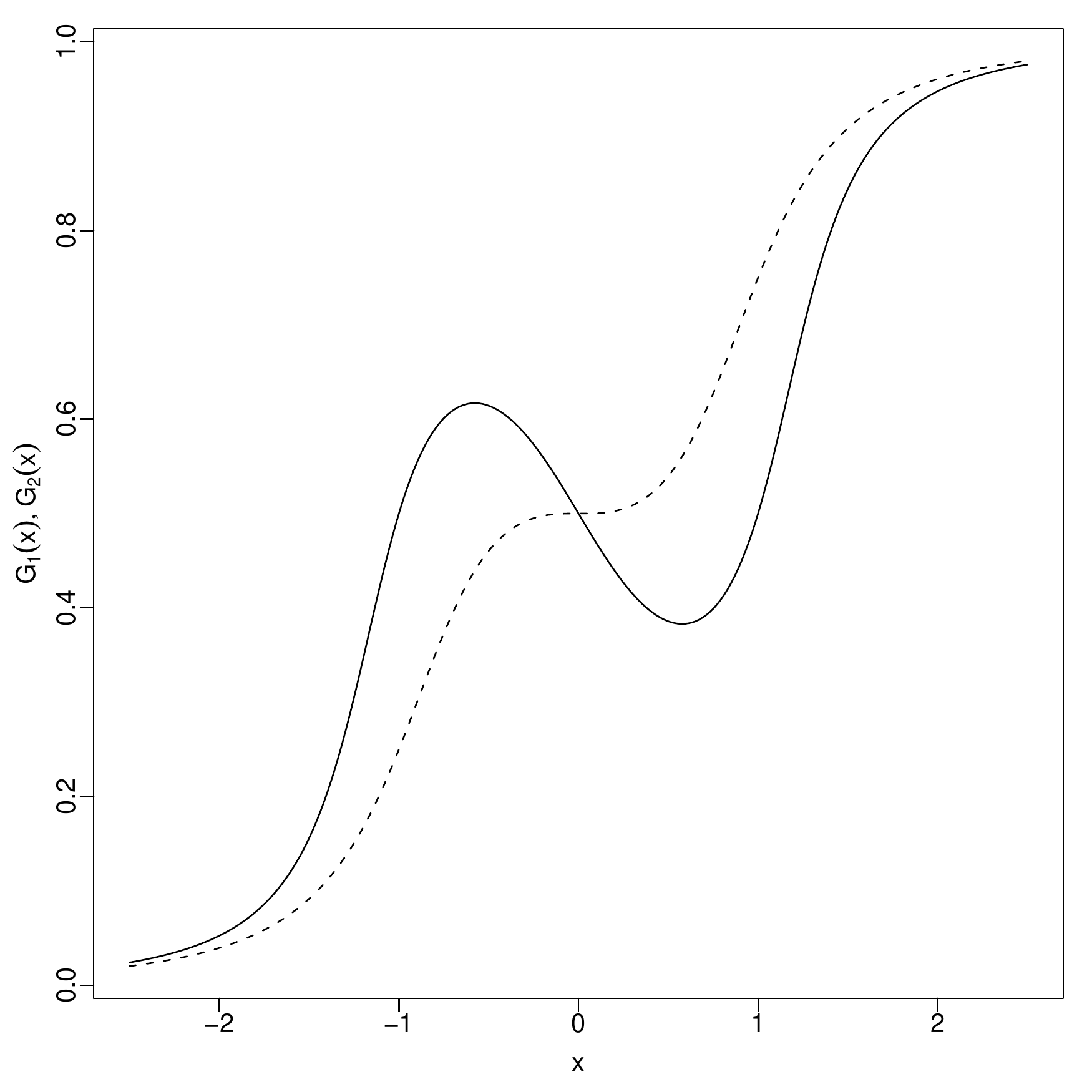} 
   \includegraphics[width=0.48\hsize]{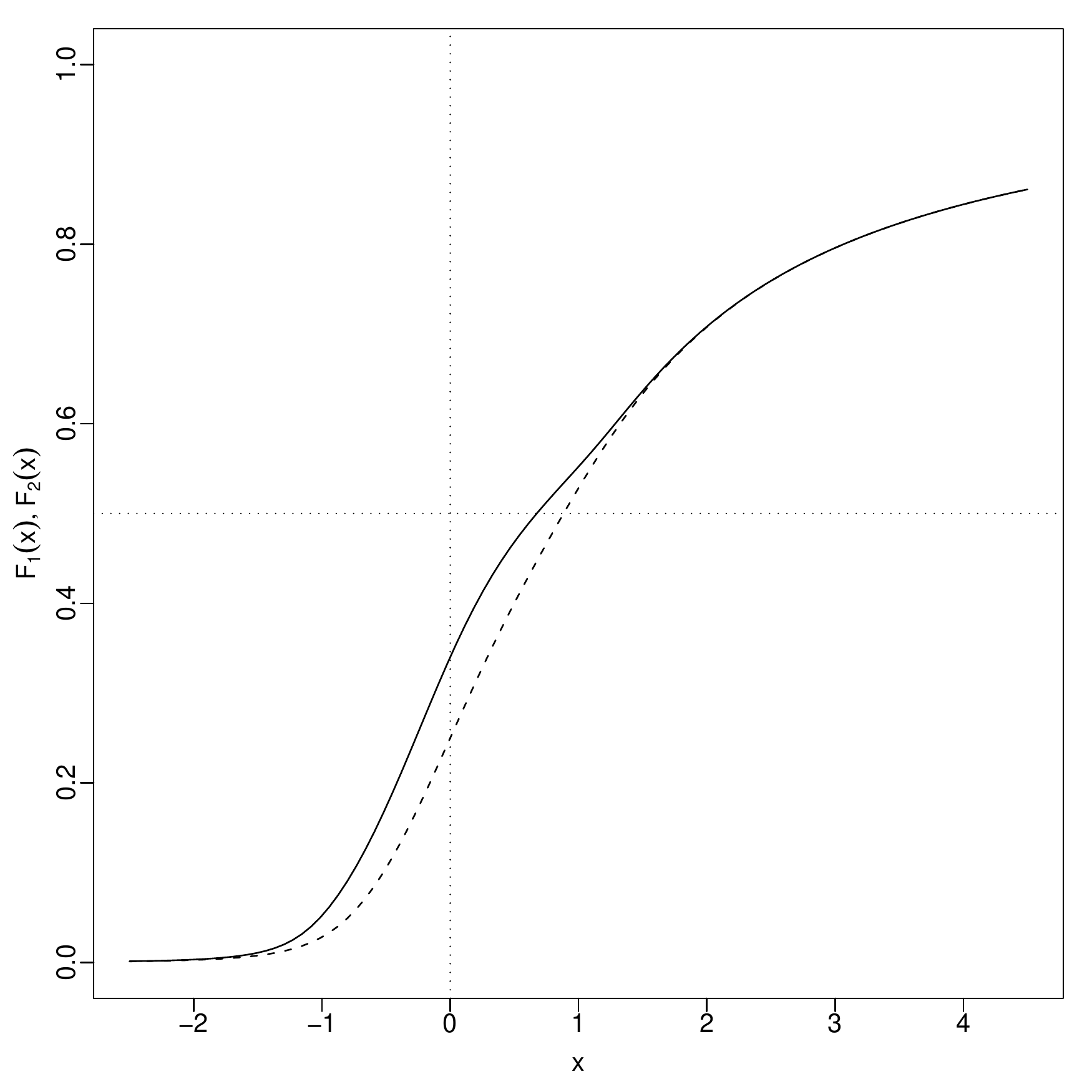} 
   \caption{Cauchy density function $f_0(x)$ perturbed by 
    $G_0$ equal to the Cauchy distribution, choosing 
    $w_1(x)=x^3-x$, and $w_2(x)=x^3$;
    on the left panel $G_1$  (continuous line) and $G_2$  (dashed line), 
    on the right $F_1$  (continuous)   and $F_2$  (dashed)
         }
   \label{fig:figure01}
\end{figure}

Figure~\ref{fig:figure01} illustrates the order $G_2 \GR G_1$ and the
stochastic order between the corresponding distributions functions
$F_1(x) \ge F_2(x) $, as stated by Proposition \ref{th:ordering}.
Here $f_0$ is the Cauchy density, $G_0$ is the Cauchy distribution
functions, and two forms of $w(x)$ are considered, namely $w_1(x)=x^3-x$,
$w_2(x)=x^3$. The two   perturbation functions $G_1(x) =G_0(w_1(x))$
and $G_2(x) =G_0(w_2(x))$ are plotted in the left panel; the right panel
displays the corresponding distribution functions $F_1(x)$ and $F_2(x)$.
 
The stochastic ordering of the $F_k$'s translates immediately into a set of
implications about ordering of moments and quantiles of the $F_k$'s. 
Specifically, if $X_k$ is a random variable with distribution function $F_k$, 
for $k=1,2$, then the following statements hold.
\begin{itemize}
\item 
  If $Q_k(p)$ denotes $p$-th quantile of $X_k$ for any $0<p<1$, 
  then 
  \[  Q_1(p) \leq Q_2(p) \,, \]
   and there exists at least one $p$ for which the inequality is strict.
\item 
  For any non-decreasing function $t$ such that the expectations exist,
  \begin{equation}  \label{e:E1<E2}
     \E{t(X_1)} \le \E{t(X_2)} 
  \end{equation}
  and the inequality is strict if $t$ is increasing. 
\end{itemize}

A further specialized case occurs when $t(x)=x^{2n-1}$ in (\ref{e:E1<E2}), for
$n=1,2,\dots$, which corresponds to the set of odd moments. In this case,
(\ref{e:E1<E2}) improves a result of \citet{umbach:2006} stating that
\[  
    \E{X_0^{2n-1}} \le \E{X_1^{2n-1}} \le \E{X_*^{2n-1}}
\]
where $X_0$ has density $f_0$ and $X_*$ has density $2f_0$ on the positive
axis, which corresponds to $G(x)=I_{[0,\infty)}(x)$ in (\ref{e:prop1-G}) and
the density of $X_1$ corresponds to a $G_1$ which is a distribution function.

It can be noticed that, if $G_2 \GR G_1 \GR G_+\equiv\half$, then
 the variances of the corresponding variables $X_k$ decrease with respect to  
 $\GR$, that is $\var{X_2}\le\var{X_1}\le \var{X_0}$, while the reverse 
holds if  $G_+\GR G_1\GR G_2$.

A simple but popular setting where Proposition~\ref{th:ordering} applies is
when $w(x)=\alpha x$, for some real $\alpha$, leading to the following
immediate implication.
\begin{proposition}
 If $f_0$ and $G_0$ are as in Lemma~\ref{th:prop1}, then the set of densities
\begin{equation} \label{e:w-linear}
   f(x;\alpha) = 2\,f_0(x)\,G_0(\alpha x) 
\end{equation}
indexed by the real parameter $\alpha$ are associated to distribution
functions which are stochastically ordered with $\alpha$.
\end{proposition}
Notice that, when $\alpha$ in (\ref{e:w-linear}) is positive, it
has a direct interpretation as an inverse scale parameter for $G_0$,
while it acts as a shape parameter for $f(x)$.
Another case of interest is given by
\[ 
   w(x) = \alpha x \sqrt{\frac{\nu+1}{\nu+x^2}}\,,
\]
which occurs in connection with the skew Student's $t$ distribution with $\nu$
degrees of freedom, studied by \citet{azza:capi:2003} and others, where $f_0$
and $G_0$ are of Student's $t$ type with $\nu$ and $\nu+1$ degrees of
freedom, respectively.
Because of Proposition~\ref{th:ordering-w}, the distribution functions 
associated to (\ref{e:prop1}) with this choice of $w(x)$ are 
stochastically ordered with respect  to $\alpha$, whether or not $f_0$
and $G_0$ correspond to a Student's $t$ distribution.

\subsection{On the uniqueness of the mode}
To examine the problem of the uniqueness of the mode of $f$
when $d=1$, it is equivalent and more convenient to study $\log f$.
If $f_0'(x)$ and $g(x)=G'(x)$ exist, then
\begin{eqnarray*}
 h(x) &=& \frac{\d}{\d x} \log f(x)\\
      &=& \frac{f_0'(x)}{f_0(x)}+ \frac{g(x)}{G(x)} \\
      &=& -h_0(x) + h_g(x)\,,
\end{eqnarray*}
say. The modes of $f$ are a subset of the solutions of the equation
\begin{equation}
    h_0(x) = h_g(x)  \,,\label{e:mode}
\end{equation}
or they are on the extremes of the support, if it is bounded. 
Since at least one mode always exists, we look for conditions to rule out 
the existence of additional modes. 

For the rest of this subsection, we assume that $G(x)$ is a monotone 
function satisfying (\ref{e:G}). Without loss of generality,
we  deal with the case that $G$ is monotonically increasing; 
for decreasing functions, dual conclusions hold.

In most common cases, $f_0$ is unimodal at 0,  hence non-decreasing
for $x\leq0$. Therefore the product $f_0(x)\,G(x)$ is increasing, and
no negative mode can exist.  The same conclusion holds if $f_0$ is 
increasing  and $G(x)$ is non-decreasing for $x\leq0$.

To ensure that there is at most one positive mode, 
some additional conditions are required. For simplicity of argument, 
we assume that  $f_0$ and $G$ have continuous derivative everywhere on the 
support $S_0$ of $f_0$; this means that we are concerned with uniqueness of
the solution of (\ref{e:mode}).
A sufficient set of conditions for this uniqueness 
is that $h_0(x)$ is increasing and $g(x)$ is decreasing.
These requirements imply that $h_g(0) >0$ and $h_g$ is decreasing, 
so that $0=h_0(0)< h_g(0)$ and the two functions can cross at most once
for $x>0$.
When  $S_0$ is unbounded, a solution of (\ref{e:mode}) always exists,
since $g\to0$ and $h_g\to 0$ as $x\to\infty$.
If $S_0$ is bounded,  (\ref{e:mode}) may happen to have no solution;
in this case, $f(x)$ is increasing for all $x$ and its mode occurs 
at the supremum of $S_0$. 
We summarize this discussion in the following statement.

\begin{proposition} \label{th:mode} If $G(x)$ in (\ref{e:prop1-G}) is a
  increasing function and $f_0(x)$ is unimodal at $0$, then no negative mode
  exists.  If we assume that $f_0$ and $G$ have continuous derivative
  everywhere on the support  $f_0$, $G(x)$ is concave for $x>0$, and
  $f_0(x)$ is log-concave, where at least one of these properties holds in a
  strict sense, then there is a unique positive mode of $f(x)$.  If $G(x)$ is
  decreasing, similar statements hold with reversed sign of the mode;
  uniqueness of the negative mode requires that $G(x)$ is convex for $x<0$.
\end{proposition}

Recall that the property of log-concavity of a univariate density function 
is equivalent to strong unimodality; see for instance Section~1.4 of 
\citet{dhar:joag:1988}. 

To check the above conditions in specific instances, it is convenient to work
with the functions $h_0$ and $g'$, if the latter exists. 
In the case of increasing $w(x)$, uniqueness of the mode is ensured 
if $g'(x)<0$ for $x>0$ and $h_0(x)$ is an increasing positive function. 
In the linear case $w(x)=\alpha\,x$, log-concavity of $f_0(x)$ and unimodality
of $g_0(x)$ at $0$ suffice to ensure unimodality of $f(x)$.
 
Table~\ref{t:f0,h0} recalls some of the more commonly employed density
functions $f_0$ and their associated functions $h_0$ and $h_0'$.
\begin{table}
\caption{Some commonly used densities $f_0$ and associated components}
\label{t:f0,h0}
\begin{center}
\begin{tabular}{lccc}
  \hline
  distribution & $f_0(x)$ & $h_0(x)$ & $h_0'(x)$ \\
  \hline
 standard normal 
    & $\phi(x)$ & $x$  &  $1$\\[1ex]
 logistic  
    & $\dfrac{e^x}{(1+e^x)^2}$ & $\dfrac{e^x -1}{e^x+1}$ 
    & $\dfrac{2\,e^x}{(e^x+1)^2}$ \\[2ex]
 Subbotin
    & $c_\nu\exp\left(-\dfrac{|x|^\nu}{\nu}\right)$ & $\sign(x)\,|x|^{\nu-1}$ 
    & $\sign(x)\,(\nu-1)|x|^{\nu-2}$ \\[2ex]
  Student's $t_\nu$ 
    & $ c_\nu \left(1 + \dfrac{x^2}{\nu}\right)^{-\frac{\nu+1}{2}}$
    & $\dfrac{\nu+1}{\nu} \:\dfrac{x}{1+x^2/\nu}$ 
    & $\dfrac{(\nu+1)(\nu-x^2)}{(\nu+x^2)^2}$ \\[2ex]
  \hline
\end{tabular}
\end{center}
\end{table}
For the first two distributions of Table~\ref{t:f0,h0},
and for the Subbotin's distribution when $\nu>1$, $h_0$ is increasing. 
If one combines one of these three choices of $f_0$ with  the distribution
function of a symmetric density  having unique mode at $0$, then 
uniqueness of the mode of $f(x)$ follows. 
Clearly, the condition of unimodality of $g(x)$ holds 
if $g_0$ is unimodal at 0 and  $w(x)=\alpha x$. 
The criterion of Proposition~\ref{th:mode} does not apply for the Student's
distribution, since $h_0(x)$ is increasing only in the interval 
$(-\sqrt{\nu}, \sqrt{\nu})$. Hence a second intersection with $h_g$ cannot be
ruled out even if $g(x)$ is decreasing for all $x>0$.
However, for the  skew-$t$ distribution, unimodality has been
established in the multivariate case by \citet{capitanio:2008} and 
\citet{jama:bala:2010jmva}, 
and furthermore it follows as a corollary of a stronger result to be 
presented in Section~\ref{s:d-general}.

The requirement of differentiability of $f_0$ and $G$ in
Proposition~\ref{th:mode} rules out a limited number of practically relevant
cases. For this reason, we did not dwell on a specific discussion of less
regular cases. One of the very few relevant distributions which are excluded 
occurs when $f_0$ is the Laplace density function.  This case is however
included in the discussion of the multivariate Subbotin distribution,
developed in Section~\ref{s:log-conc}, when $\nu=1$ and $d=1$.

\begin{figure}
  \centering
   \includegraphics[width=0.48\hsize]{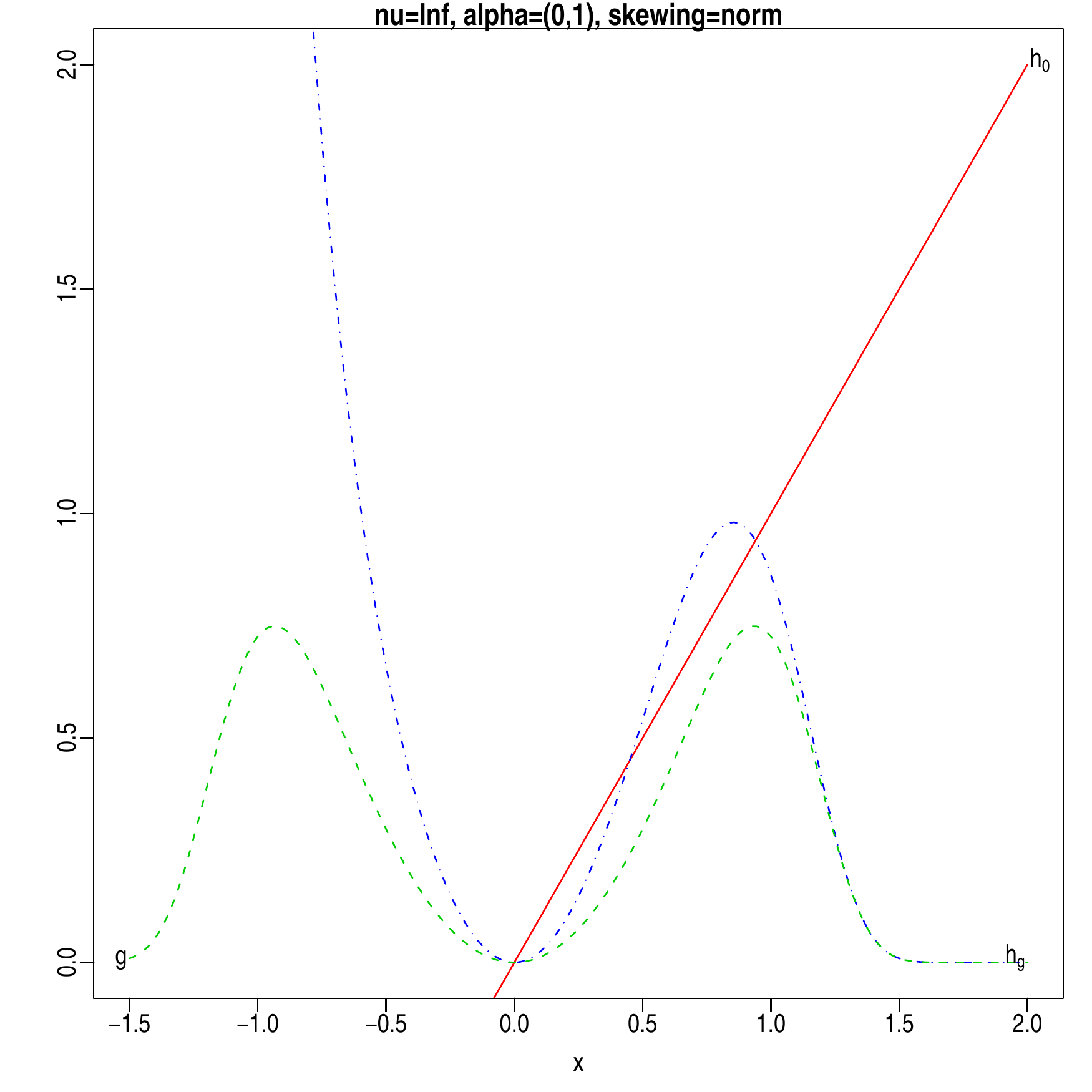} 
   \includegraphics[width=0.48\hsize]{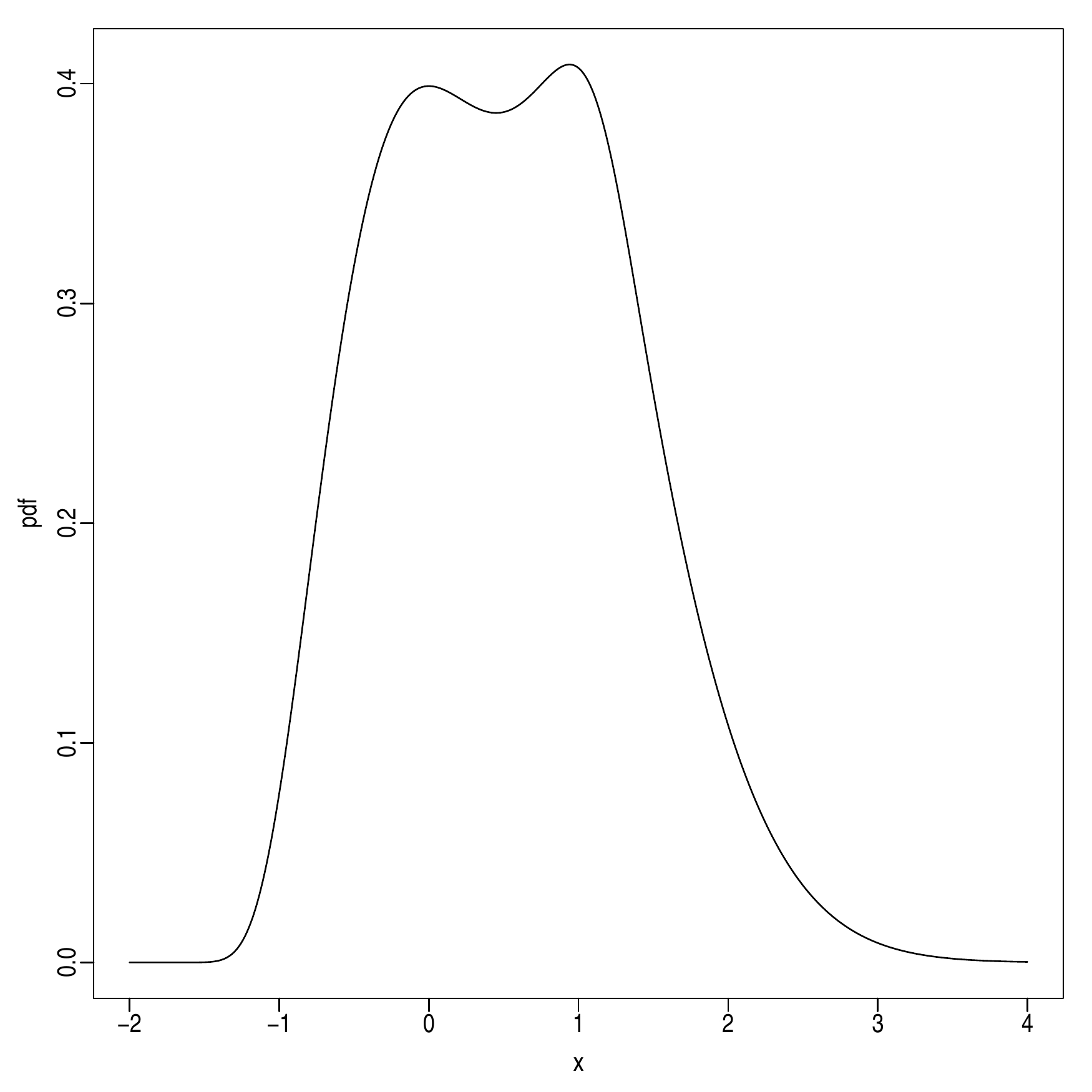} 
   \caption{Case with $f_0=\phi$, $G_0=\Phi$ , $w(x)=x^3$.  The left-side
     panel displays the functions $h_0$ (continuous), $g$ (dashed), $h_g$
     (dot-dashed); the right-side panel displays $f(x)$}
   \label{f:multimode}
\end{figure}

Although Proposition~\ref{th:mode} only gives a set of sufficient conditions
for unimodality, the condition that $g(x)$ is decreasing for $x>0$ cannot be
avoided completely. In other words, when $f$ is represented in the form
(\ref{e:prop1}), the sole condition of increasing $w(x)$ is not sufficient for
unimodality. This fact is demonstrated by the simple case with $f_0=\phi$,
$G_0=\Phi$, $w(x)=x^3$, whose key features are illustrated in
Figure~\ref{f:multimode}. Since $w'(0)=0$, then $g(0)=0$; hence (\ref{e:mode})
has a solution in 0, but the left panel of Figure~\ref{f:multimode} shows that
there are two more intersections of $h_0$ and $h_g$ for $x>0$, one
corresponding to an anti-mode and one to a second mode of $f(x)$, as visible
from the right panel of the figure.

This case falls under the setting examined by \citet{maY:gent:2004}
who have shown that for $f_0=\phi$, $G_0=\Phi$, 
$w(x)=\alpha\,x +\beta\,x^3$ there are at most two modes.
Some additional  conditions may ensure unimodality: one such set
of conditions is $\alpha,\beta>0$ and $\alpha^3>6\beta$. To prove
that they imply unimodality of $f$, consider
\begin{eqnarray*}
\frac{\d^2 \log\Phi(w(x))}{\d x^2} &=&
   -\frac{\phi(w(x))}{\Phi(w(x))^{2}} \\
&& \times
 \left\{\Phi(w(x))[(\beta x^3+\alpha x)(3\beta x^2+\alpha)^2-6\beta\,x]
 +\phi(w(x))(3\beta x^2+\alpha)^2\right\}.
\end{eqnarray*} 
whose terms inside curly brackets, except $-6\beta\,x$, are all positive for
$x\geq0$.  Since $\alpha^3>6\beta$, then $(\alpha^3-6\beta) x\,$ is positive,
so that this derivative is negative and $G_0(w(x))$ is log-concave for 
$x\geq 0$.  
For $x<0$, we use this other argument: since $G_0$ is increasing and
log-concave and $w(x)$ is concave in the subset $x<0$, then the composition
$G_0(w(x))$ is log-concave in the subset $x<0$; see
Proposition~\ref{th:MO}\,(iii) below.  Since the above second derivative is
continuous everywhere, then $G_0(w(x))$ is log-concave everywhere.

\section{Quasi-concave and unimodal densities in $d$ dimensions} 
\label{s:d-general}
 
A real-valued function $f$ defined on a subset $S$ of $\Real^d$ is said to be
quasi-concave if the sets of the form $C_u=\{x: f(x)\ge u\}$ are convex for
all positive $u$. If $d=1$, the notion of quasi-concavity coincides with 
uniqueness of the maximum, provided a pole is regarded as a maximum point,
but for $d>1$ the two concepts separate.
This motivates the following digression about concavity and related concepts, 
to develop some tools which will be used later on for our main target.

\subsection{Concavity, quasi-concavity and unimodality} 

We first recall some standard notions available for instance in Chapter~16
of \citet{marshall:olkin:1979}.
A real function $f$ defined on a convex subset $S$ of
$\Real^d$ is said to be concave if, for every $x$ and $y \in S$ and 
$\theta \in (0,1)$, we have
\[ 
   f(\theta x+ (1-\theta) y) \geq \theta f(x)+ (1-\theta) f(y);
\]
in this case $-f$ is a  convex function. A function $f$ is  said to be 
log-concave if $\log f$ is concave, that is for every $x$ and $y 
\in S$ and $\theta \in (0,1)$ we have
\[
   f(\theta x+ (1-\theta) y) \geq  f(x)^{\theta}  f(y)^{1-\theta}.
\]
The terms strictly concave and strictly log-concave apply if the above
inequalities  hold in a strict sense for all $x\neq y$ and all $\theta$.

Concave and log-concave functions defined on  an open set are continuous. 
Moreover a twice differentiable function is concave (strictly concave) 
if and only if its Hessian matrix is negative semi-definite 
(negative definite) everywhere on $S$.

The next proposition provides the concave and log-concave extension of 
classical composition properties for convex functions such as statement 
(i) which can be found for example in \citet[p.\,451]{marshall:olkin:1979}
together with its proof; the proofs of the
other statements are completely analogous.

\begin{proposition} \label{th:MO} Let $h$ be a real function defined on a
  convex set $S$, a subset of $\Real^d$, and $H$ a monotone real function
  defined on a convex subset of $\Real$, such that the composition
  $H(h)$ is defined on $S$.  Then the following properties hold.
  \begin{enumerate}[(i)]
  \item If $h$ is convex and $H$ non-decreasing and convex, then
    $H(h)$ is convex. Moreover $H(h)$ is strictly convex if
    $H$ is strictly convex, or if $h$ is strictly convex and
    $H$ is strictly monotone.
 
  \item If $h$ is convex and $H$ non-increasing and log-concave, then
    $H(h)$ is log-concave. Moreover $H(h)$ is strictly
    log-concave if $H$ is strictly log-concave, or if $h$ is strictly
    convex and $H$ is strictly monotone. The same statements hold
    replacing the term log-concave by concave throughout.

  \item If $h$ is concave and $H$ non-decreasing and log-concave, then
    $H(h)$ is log-concave. Moreover $H(h)$ is strictly
    log-concave if $H$ is strictly log-concave, or if $h$ is strictly
    concave and $H$ is strictly monotone. The same statements hold
    replacing the term log-concave by concave throughout.
  \end{enumerate}
\end{proposition}

We have defined quasi-concavity by requiring convexity of all
sets $C_u$. An equivalent condition is that,  for every $x$ and
$y \in S\subseteq \Real^d$ and $\theta \in (0,1)$, we have
\[
   f(\theta x+ (1-\theta) y) \geq  \min\{f(x), f(y)\}.
\]
Obviously a function which is concave or log-concave is also quasi-concave.
Similarly, both strict concavity and strict log-concavity imply strict
quasi-concavity.

We now apply the above notions to the case where $f$ represents a probability
density function on a set $S\subseteq \Real^d$.  The concept of unimodality
has a friendly formal definition in the univariate case, see for
instance \citet[p\,.2]{dhar:joag:1988}, but this has has no direct equivalent
in the multivariate case. Informally, we say that\ the term mode of a density
refers to a point  where the density takes a maximum
value, either globally or locally.  While a boring formal definition which
allows for the non-uniqueness of the density function could be given, such a
definition is not really necessary for the main aims of the present paper,
since the density functions which we are concerned with are so regular that
their modes are either points of (local) maxima or poles.

The set of the modes of a quasi-concave density is a convex set.  Moreover, 
if $f$ is strictly quasi-concave, then the mode is unique.  When the mode is
unique we say that density $f$ is \emph{unimodal}, and we say that $f$ is
\emph{c-unimodal} if the set of its modes is a convex set.
If $X$ is a random variable with density function $f$ which is
unimodal, we shall say that $X$ is unimodal, with slight abuse of terminology.
The same convention is adopted for log-concavity, quasi-concavity and other 
properties.

Another important notion is $s$-concavity, which helps to make the concept of
quasi-concavity more tractable.  A systematic discussion of $s$-concavity has
been given by \citet{dhar:joag:1988}; see specifically their Section~3.3, of
which we now recall the main ingredients.  Given a real number $s\neq0$, a
density is said to be $s$-concave on $S$ if
\[
  f(\theta x+ (1-\theta) y) \geq \{\theta f(x)^s+ (1-\theta) f(y)^s\}^{1/s}.
\]
for all $x,\ y \in S$ and  all $\theta\in(0,1)$.

Clearly, concavity corresponds to $s=1$.  A density $f$ is $s$-concave with
$s<0$ if and only if $f^s$ is convex; similarly, a density $f$ is $s$-concave
with $s>0$ if and only if $f^s$ is concave.  If we call ($-\infty$)-concave a
function which is quasi-concave and $0$-concave a function which is
log-concave, then the class of sets of $s$-concave functions is increasing
when $s$ decreases; in other words, if $f$ is $s$-concave, then it is
$r$-concave for any $r<s$.  Finally, notice that is easy to adapt Proposition
\ref{th:MO} to $s$-concave functions.

The closure with respect to marginalization of $s$-concave densities depends
on the value of $s$ and on the dimensions of the spaces, as indicated by the
next proposition, which essentially is Theorem~3.21 of \citet{dhar:joag:1988}.

\begin{proposition} \label{th:marg} 
  Let $f$ be an $s$-concave density on a convex set $S$ in $\Real^{d+m}$,
  and $f_d$ be the marginal density of $f$ on an $d$-dimensional subspace.
  If $s\geq -1/m$, then $f_d$ is $s_m$-concave
  on the projection of the support of $f$, where $s_m=s/(1+ms)$, with the
  convention that, if $s=-1/m$, then $s_m=-\infty$.
\end{proposition}
Notice that this result includes the fact that the class of log-concave 
densities is closed with respect to marginalization. 
In addition,  from a perusal of the proof of the above-quoted Theorem~3.21,
we obtain that the marginal densities are strictly $s_m$-concave provided 
$f$ is strictly $s$-concave or the set $S$ is strictly convex.


\subsection{Skew-elliptical distributions generated by conditioning}
\label{s:skew-ell}
 
A $d$-dimensional random variable $U$ is said to have an elliptical density, 
with density generator function $\tilde f$, if its density $f_U$ is of the form
\begin{equation} \label{e:SE}
  f_U(y)= k\, \tilde f(y\T\Omega^{-1}y),
\end{equation}
where $\Omega$ is a $d$-dimensional positive definite  matrix, 
the function $ \tilde f:(0,+\infty)\to \Real^+$ is such that 
$x^{d/2-1} \tilde f(x)$ has finite integral on $(0,+\infty)$ and 
$k$ is a suitable constant which depends on $d$ and $\det(\Omega)$. 
In this case, we shall use the notation $U\sim\mathcal{E}_d(0,\Omega,\tilde f)$.

Note that an elliptical density $f$ is c-unimodal if and only if its density
generator is non-increasing, and it is unimodal if and only if its density
generator is decreasing. 
Then it turn out that $f$ is c-unimodal if and only if it is quasi-concave, 
and it is unimodal if and only if it is strictly quasi-concave.

An initial formulation of skew-elliptical distribution has been considered
by \citet{azza:capi:1999}, which was of type (\ref{e:prop1}) with $f_0$ of
elliptical class and $w(x)$ linear.  Another formulation 
of skew-elliptical distribution has been put forward by \citet{bran:dey:2001},
whose key ingredients are now recalled.
Consider a $(d+1)$-dimensional random variable
\begin{equation} \label{e:U}
  U=\pmatrix{U_0\cr U_1}\sim\mathcal{E}_{d+1}(0,\Omega_+, \tilde{f})\,,
 \qquad \mathrm{where}\quad
 \Omega_+=\pmatrix{1&\delta\T\cr\delta& \Omega} > 0,
\end{equation}  
and $U_0$ and $U_1$ have dimension $1$ and  $d$, respectively; for our
aims, there is no loss of generality in assuming that the diagonal elements 
of $\Omega_+$ are all $1$'s. Then a random variable  $Z= (U_1|U_0>0)$ is 
said to have a  skew-elliptical distribution,  and its density function at 
$u_1\in\Real^d$ is
\begin{equation} \label{e:SE-pdf}
  f_Z (u_1)=2\int_0 ^{+\infty} k_1 \tilde{f}(u\T\Omega_+\inv{}u)\d u_0
\end{equation}
where $u\T=(u_0,u_1\T)$. This construction arises as an extension 
of one of the mechanisms for generating the skew-normal distribution
to the case of elliptical densities, but the study of the connections 
with other densities of type (\ref{e:prop1}) was not an aim of
\citet{bran:dey:2001}.

Consequently, one question investigated by \citet{azza:capi:2003} 
was whether all distributions
of type (\ref{e:SE-pdf}) are of  type (\ref{e:prop1}),  
with the requirement that $f_0$ is the density of an elliptical 
$d$-dimensional distribution.
The conjecture has been proved for a set of important cases, notably the
multivariate skew-normal and the skew-$t$ distributions, among others, but
a general statement could not be obtained.
This general conclusion is however quite simple to reach using 
representation (\ref{e:prop1-G}), and recalling that 
\citet{bran:dey:2001} have proved that  (\ref{e:SE-pdf}) can be written as
\begin{equation} \label{e:branco-dey}
    f_Z (y) = 2\,\ f_0(y)\,F_y(\alpha\T\,y),  \qquad\qquad (y\in\Real^d),
\end{equation}
where $f_0$ is the density of an elliptical $d$-dimensional distribution, and
$F_y$ is a cumulative distribution function of 
a symmetric univariate distribution, which
depends on $y$ only through $y\T \Omega^{-1}y$. Since $F_y = F_{-y}$, 
then it is immediate that $G(y)= F_y(\alpha\T\,y)$ satisfies (\ref{e:G}).
Hence (\ref{e:branco-dey}) allows a representation of type (\ref{e:prop1-G}),
and via (\ref{e:G0-G}) also of type (\ref{e:prop1}).

\begin{proposition} \label{th:propSE} Assume that the random variable $U$ in
  (\ref{e:U}) is c-unimodal.  If $\tilde{f}$ is log-concave, then 
  the  elliptical densities of $U$ and $U_1$ and the skew-elliptical density 
  of $Z$ are log-concave.  
  Moreover they are strictly log-concave if $U$ is unimodal or $\tilde{f}$ is
  strictly log-concave or the support of $\tilde{f}$ is   bounded.
\end{proposition} 

\noindent\textbf{Proof}. 
Function $h(u)= u\T\Omega_+\inv{}u$ is strictly convex.  Since $U$ is
c-unimodal, then $\tilde{f}$ is non-increasing, moreover it is log-concave;
therefore $\tilde{f}(u\T\Omega_+\inv{}u)$ is log-concave by
Proposition~\ref{th:MO}\,(ii).  Then both $U$ and $(U|U_0>0)$ have log-concave
densities.  Since the marginals of a log-concave density are log-concave, then
log-concavity of $U_1$ and $Z$ holds by (\ref{e:SE-pdf}).  Now, if $U$ is
unimodal, $\tilde{f}$ is decreasing, and $\tilde{f}(u\T \Omega_+\inv{}u)$ is
strictly log-concave, by Proposition~\ref{th:MO} (i).  If $\tilde{f}$ is
strictly log-concave, then $\tilde{f}(u\T\Omega_+\inv{}u)$ is strictly
log-concave. Finally, if the support of $\tilde{f}$ is bounded, then the
support of $U$ is strictly convex and, by Proposition~\ref{th:MO} (i), also in
this case $\tilde{f}(u\T \Omega_+\inv{}u)$ is strictly log-concave.  Then, in
all three cases, strict log-concavity of $U_1$ and $Z$ holds by recalling the
remark following Proposition~\ref{th:marg}.

This proposition is a special case of the more general result which follows,
but we keep Proposition~\ref{th:propSE} separate both because of the special
role of log-concavity and because this arrangement allows a more compact
exposition of the combined discussion.

\begin{proposition} \label{th:propSEsconc} Assume that the random variable $U$
  in (\ref{e:U}) is c-unimodal.  If $\tilde{f}$ is $s$-concave, with 
  $s\geq -1$, then    $U$  has s-concave density, whereas
  the elliptical density of $U_1$ and the skew-elliptical density of
  $Z$ are $s_1$-concave, with $s_1=s/(1+s)$.  
  Moreover  all conclusions hold  strictly  if $U$ is unimodal or $\tilde{f}$ 
  is strictly $s$-concave or the support of $\tilde{f}$ is bounded.
 \end{proposition} 

\noindent\textbf{Proof}.  
The function $h (u)= u\T\Omega_+\inv{}u$ is strictly convex.
Since $U$ is c-unimodal, then $\tilde{f}$ is non-increasing and moreover 
it is $s$-concave. 
We now examine properties of concavity  separating the case $s<0$
and $s>0$; the case $s=0$, which corresponds to log-concavity, has already 
been handled in Proposition~\ref{th:propSE}.
If $s<0$ then $\tilde{f}^s$ is non-decreasing and convex.  
Then $\tilde{f}^{s}(u\T\Omega_+\inv{}u)=\{\tilde{f}(u\T\Omega_+\inv{}u)\}^s$
is convex by Proposition~\ref{th:MO}\,(i)
and $\tilde{f}(u\T\Omega_+\inv{}u)$ is $s$-concave.  
On the other hand, if $s>0$ then $\tilde{f}^s$ is non-increasing and concave.  
Then $\tilde{f}^{s}(u\T\Omega_+\inv{}u)=\{ \tilde{f}(u\T\Omega_+\inv{}u)\}^s$ 
is concave by Proposition~\ref{th:MO}\,(ii) and $\tilde{f}(u\T\Omega_+\inv{}u)$
is $s$-concave.  Then both $U$ and $(U|U_0>0)$ have $s$-concave densities.
Now, the claim about the densities  of $U_1$ and  $Z$ follows from  
Proposition \ref{th:marg} by taking into account (\ref{e:SE-pdf}).  
The final statement follows by the same type of argument used in the proof 
of Proposition~\ref{th:propSE}. 

Note that, in the special case of a concave density generator, the support 
is bounded, and both the marginal density on $\Real^d$ and the skew-symmetric 
density of $Z$ are not necessarily concave. However, using
Proposition~\ref{th:propSEsconc} with $s=1$,  strict $1/2$-concavity of their
densities follows, and this fact implies strictly log-concavity. 

The results of Proposition~\ref{th:propSE} and
Proposition~\ref{th:propSEsconc} allow to handle several classes of
distributions, of which we now sketch the more noteworthy cases.

A important specific instance is the multivariate skew-normal density which 
can be  represented by a conditioning method.
For an expression of the multivariate skew-normal density, 
see for instance (16) of \citet{azzalini:2005}.
 Since the density generator of
the normal family, $\tilde{f}(x)= \exp (-x/2)$, is decreasing and log-concave, 
then from Proposition~\ref{th:propSE} we obtain log-concavity of the 
skew-normal family. This conclusion is however a special case of a more
general result on log-concavity of the SUN distribution obtained by
\citet{jama:bala:2010jmva}; see their Theorem~1.

The $(d+1)$-dimensional Pearson type II distributions for which 
$\tilde{f}(x)= (1-x)^{\nu}$, where $x\in (0,1)$ and  $\nu\geq 0$, 
satisfies the conditions of Proposition \ref{th:propSEsconc}. 
In fact it is non-increasing and $\nu\inv$-concave on a bounded support. 
Then the skew-elliptical $d$-dimensional density is strictly 
$(\nu+1)\inv$-concave and therefore strictly log-concave.  
 The density function of the skew-type\,II density function
is given by (22) of \citet{azza:capi:2003}.

In addition, Proposition \ref{th:propSEsconc} holds for the Pearson type~VII 
distributions, and in particular for the Student's distribution. 
In this case  the density generator is given by
\begin{equation}
     \tilde f(x)=(1+x/\nu)^{-M} \qquad\qquad (\nu<0,\; (d+1)/2<M).
\end{equation}
and $M=(d+\nu+1)/2$ for the Student's density.  Such generator is decreasing
and $s$-concave with $s=-1/M$; in fact $\tilde f(x)^{-1/M}$ is convex.  Since
$s\geq-1$, then Proposition~\ref{th:propSEsconc} applies and the skew-$t$ is
$s_1$-concave with $s_1= -1/(M-1)$, and $s_1= -2/(d+\nu-1)$ in the Student's
case. These densities are not log-concave, but they are still strictly
quasi-concave. Hence unimodality follows.   For 
expressions of the multivariate skew-type\,VII and skew-$t$ density, 
see (21) and (26) of \citet{azza:capi:2003}, respectively.

The above results establish not only unimodality of the more appealing subset
of the skew-elliptical family of distributions, namely those of type
(\ref{e:SE-pdf}), but also the much stronger conclusion of quasi-concavity of
these densities.  It is intrinsic to the nature of skew-elliptical densities
that they do not have highest density regions of elliptical shape, but it is
reassuring that they maintain a qualitatively similar behaviour, in the sense
that convexity of these regions, $C_u$ in our notation, holds as long as the
parent $(d+1)$-dimensional elliptical density enjoys a qualitatively similar
property but in a somewhat stronger variant, specifically $s$-concavity
with $s\ge-1$.

Note that there is no hope to extend Proposition \ref{th:propSE} to 
quasi-concave densities,  in the sense that a skew-symmetric generated by 
conditioning  a quasi-concave density is not  necessarily quasi-concave as
demonstrated by the following construction.

\begin{figure}
  \centering
   \includegraphics[width=0.48\hsize]{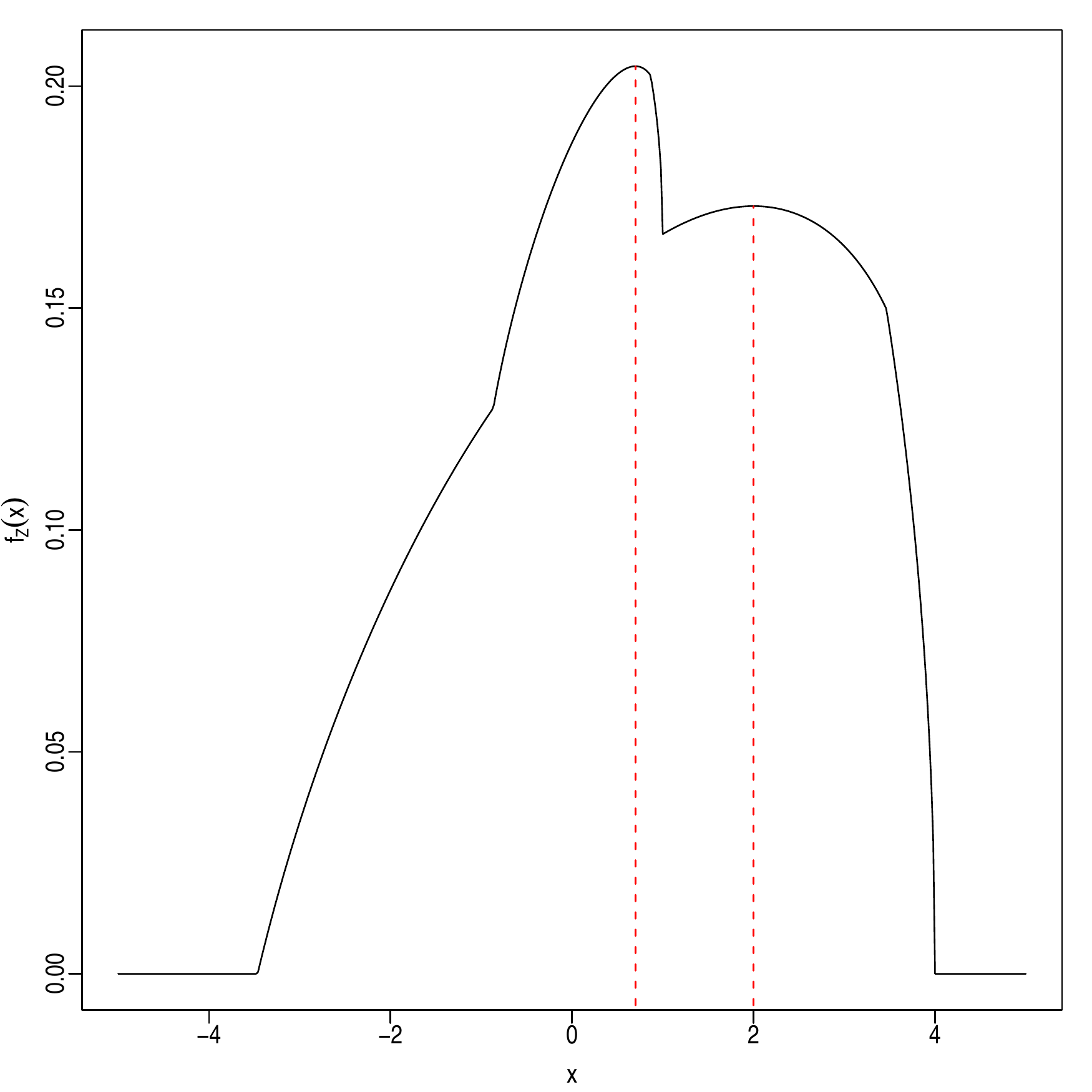}
   \caption{Density function $f_Z(x)$, exhibiting lack of quasi-concavity,
     obtained by conditioning of a bivariate elliptical quasi-concave 
     distribution}
   \label{f:lalinea}
\end{figure}

\paragraph{Example} 
Consider $U= (U_0, U_1)\T \sim\mathcal{E}_2(0,\Omega_+,\tilde f)$, where 
\[ 
\tilde f = I_{(0,1)}+ I_{(0,4^2)}
 \qquad \mathrm{and}\quad
 \Omega_+=\pmatrix{1& 1/2\cr 1/2& 1} ,
\] whose density function is
 $$f_U(x,y)= k\{I_{S_1}(x,y)+ I_{S_4}(x,y) \}$$
where $S_j= \{(x,y)\in \Real^2 : x^2+y^2-xy\leq  3j^2/4  \}$, $j=1, 4$,
and $k$ is the normalizing constant given by 
$k = 1/(A_1 + A_4)\approx 0.0216$ where $A_j=\pi\sqrt{3}j^2/2$ . 
Then both $U_0$, and $U_1$ have common support  $[-4, 4]$ and  
density function
$$  
  f_{U_0}(x)=f_{U_1}(x) = k\, \left(
     \sqrt{3(1-x^2)}\, I_{(-1, 1)}(x)+\sqrt{3(16-y^2)} \right)\,.
$$ 
Because of (\ref{e:SE}) and (\ref{e:SE-pdf}), the density of $Z=(U_1| U_0>0)$ 
is given by
 $$f_Z(y)= k \{f_1(y)+f_4(y) \} $$
where 
\[ 
   f_j(y) = 2\int_0^{+\infty} { I_{S_j}(x,y)}\d x =  
   \cases{y+\sqrt{3(j^2-y^2)}   \quad & 
                        if $-\sqrt3 j/2 \leq y \leq \sqrt3 j/2 $, \cr
          2\sqrt{3(j^2-y^2)}     & if $\sqrt3 j /2 \leq y \leq j$, \cr
                       0        & otherwise,} \\
\]      
for $j=1,4$, and it is displayed in Figure~\ref{f:lalinea}.  
The global maximum of $f_Z$ is where 
$k(2 y+\sqrt{3(16-y^2)}+\sqrt{3(1-x^2)})$ takes its maximum value, that is at
$y\approx 0.699$.  When $y>1$, $f_Z=kf_4$ and there is another local maximum
at $y=2$. Therefore, $f_Z$ is not unimodal. 
 
To conclude with, while the density of $U$ is quasi-concave, the
skew-elliptical variable $Z$ generated by conditioning is not quasi-concave.

\subsection{Log-concavity of other families of distributions}
\label{s:log-conc}

There are several other families of distributions which belong to the 
area of interest of the stream of literature described at the beginning 
of this paper but are not included in the conditioning mechanism of
an elliptical distribution considered in \S\,\ref{s:skew-ell}.
This section deals with log-concavity of some of these other families,
making use of the following immediate implication of Proposition~\ref{th:MO}.

\begin{corollary} \label{th:log-concave} If $q_0$ is a log-concave function
  defined on a convex set $S\subseteq\Real^d$, and $H$ and $h$ are as in
  Proposition~\ref{th:MO}, either (ii) or (iii), then
\begin{equation}
  q(x) = q_0(x)\,H\{h(x)\}, \qquad (x\in S),
\end{equation}
is log-concave on $S$.
\end{corollary}

\paragraph{Example} The density function on the real line introduced by \cite{subbotin:1923} has
been variously denoted  by subsequent authors as  exponential power 
distribution, generalized error  distribution and  normal distribution 
of order $\nu$.  Its multivariate version is
\[  
     f_\nu(x) =  c_\nu\det(C)^{1/2} \,
            \exp\left(-\frac{(x\T C x)^{\nu/2}}{\nu}\right), \qquad
             (x\in\Real^d),
\] 
where $C$ is a symmetric positive definite matrix, 
$\nu$ is a positive parameter and $c_\nu$ a normalization constant.
For $\nu=2$ and $\nu=1$, $f_\nu$ lends the multivariate normal and the
multivariate Laplace density, respectively.

We first want to show that $f_\nu$ is log-concave if $\nu\ge 1$. Consider 
$h(x)=(x\T C x)^{1/2}$ whose Hessian matrix is
\[  
  \pd{^2\,h(x)}{x\:\partial x\T}= 
   h(x)^{-3}\left( x\T C x C - C x x\T C \right) =  h(x)^{-3} M,
\]
say. To show that this Hessian is positive  semi-definite,
it is sufficient to prove this fact for matrix $M$, since $h(x)\ge0$.
For any  $u\in\Real^d$, write
\[
   u\T M u 
      =  (x\T Cx) (u\T C u) -(u\T C x) (x\T C u) 
      =  \|\tilde{u}\|^2 \,\|\tilde{x}\|^2 -(\tilde{u}\T\tilde{x})^2 
\]
where $\tilde{u}=C^{1/2}u$ and $\tilde{x}=C^{1/2}x$ for any square
root $C^{1/2}$ of $C$,
and from the Cauchy-Schwarz inequality we conclude that $u\T M u \ge 0$.
Then $h$ is convex.   Next, write
\[   -\log f_\nu(x)= \mbox{constant}+h(x)^\nu/\nu \]
and observe  that,  since $t^\nu$ is a strictly convex for $t\ge0$, 
then $-\log f_\nu$ is convex for $\nu\ge 1$ and strictly convex for $\nu>1$ 
by Proposition~\ref{th:MO}\,(i). Hence $f_\nu$ is log-concave for $\nu\ge1$ 
and strictly log-concave for $\nu>1$.

Now we introduce a skewed version of $f_\nu$ of type (\ref{e:prop1}). 
If we aim at obtaining a density which fulfils the requirements
of both Lemma~\ref{th:prop1} and Corollary~\ref{th:log-concave},
then $H=G_0$  is non-decreasing, while function $h = w$ must be odd 
and concave,  hence it  has to be linear.   
We then focus on the density function 
\begin{equation} \label{e:SEP}
    f(x) = 2 \,f_\nu(x) \,G_0(\alpha\T x), \qquad (x\in \Real^d), 
\end{equation}
where $G_0$ is a distribution function on $\Real$, symmetric about $0$.

Among the many options for $G_0$, a quite natural choice is to take
$G_0$ equal to the distribution function of $f_\nu$ in the scalar case,
that is
\[
  G_0(t) = \frac{1}{2}\left(1+
        \sign(t) \frac{\gamma(|t|^\nu/\nu, 1/\nu)}{\Gamma(1/\nu)}\right) \,,
        \qquad t\in\Real,
\]
where $\gamma$ denotes the lower incomplete gamma function.
This choice of $G_0$ has been examined by \citet{azzalini:1986} in
the case $d=1$ of (\ref{e:SEP}). 
He has shown that $G_0$ is strictly log-concave if $\nu>1$, 
leading to log-concavity of  (\ref{e:SEP}) when $d=1$.
The case $\nu=1$ which corresponds to the Laplace distribution function 
is easily handled by direct computation of the second derivative to
show strict log-concavity of $G_0$.
Now, combining strict log-concavity of $G_0$ with log-concavity of $f_\nu$
proved above, an application of Corollary~\ref{th:log-concave} shows that
(\ref{e:SEP}) is strictly log-concave on $\Real^d$ if $\nu\ge 1$.

Although (\ref{e:SEP}) is of skew-elliptical type, it is not of the type
generated by the conditioning mechanism of a $(d+1)$-dimensional
elliptical variate considered in Section~\ref{s:skew-ell}. 
In fact, the results of \citet{kano:1994} show that the set of densities
$f_\nu $ is not closed under marginalization, and this fact affects the
conditioning mechanism (\ref{e:SE-pdf}) as well. 

As an example of non-elliptical  distribution, we can consider a
$d$-fold product of univariate Subbotin's densities, that is
\[ 
   f_\nu^*(x) = \prod_{j=1}^d c_\nu \exp(-|x_j|^\nu/\nu)   \,,
   \qquad x=(x_1,\dots, x_d)\in\Real^d \,,
\]
and this density can be used as a replacement of $f_\nu$ in (\ref{e:SEP}).
 Since each factor of this product is log-concave, if $\nu\ge1$,
the same property holds for $f^*_\nu$. 
Strict log-concavity holds for $2f^*_\nu(x)
G_0(\alpha\T x)$ as well, using again strict log-concavity of $G_0$.

\paragraph{Example} 
To illustrate the applicability of Corollary~\ref{th:log-concave} to 
distributions outside the set of type  (\ref{e:prop1}), 
consider the  so-called extended skew-normal density 
which in the $d$-dimensional  case takes the form  
\begin{equation} \label{e:ESN}
    f(x) = \phi_d(x;\Omega)\,\frac{ \Phi(\alpha_0 + \alpha\T\,x)}{\Phi(\tau)},
     \qquad\qquad 
    (x\in\Real^d),
\end{equation}
where $\tau\in\Real$ and $\alpha_0= \tau(\alpha\T \Omega\alpha)^{1/2}$.
Although this distribution does not quite fall under the umbrella of 
Lemma~\ref{th:prop1} unless $\tau=0$, its constructive argument is closely 
related.
 
To show log-concavity of (\ref{e:ESN}), first recall the well-known fact 
that  $\phi_d(x;\Omega)$ is strictly 
log-concave. Moreover $\Phi$ is log-concave, as if follows by direct 
calculation of the second derivative of $\log\Phi$, taking into account 
the well-known  fact $-y\Phi(y)< \phi(y)$  for every $y \leq0$.
In addition, since $\Phi$ is strictly increasing  and $\alpha_0+\alpha\T\,x$ 
is concave in a non-strict sense, Corollary~\ref{th:log-concave} applies to 
conclude that (\ref{e:ESN}) is strictly log-concave.

Although this conclusion is a special case of the result of 
\citet{jama:bala:2010jmva} concerning log-concavity of the SUN distribution, 
it has however been presented because the above argument is different.

%

\paragraph*{Acknowledgements}
This research has been supported by MIUR, Italy, under grant scheme
PRIN, project No.~2006132978.

\end{document}